\documentclass[12pt]{amsart}

\usepackage{latexsym}
\usepackage[pdftex]{graphicx}

\newtheorem{theorem}{Theorem}
\newtheorem{lemma}{Lemma}
\newtheorem{proposition}{Proposition}
\newtheorem{remark}{Remark}

\newtheorem{corollary}{Corollary}

\def\R{{\mathbb R}}

\newcommand{\beq}{\begin{equation}}
\newcommand{\eeq}{\end{equation}}
\newcommand{\beqna}{\begin{eqnarray*}}
\newcommand{\eeqna}{\end{eqnarray*}}
\newcommand{\beqn}{\begin{equation*}}
\newcommand{\eeqn}{\end{equation*}}
\newcommand{\bp}{\begin{proof}}
\newcommand{\ep}{\end{proof}}
\newcommand{\bprop}{\begin{proposition}}
\newcommand{\eprop}{\end{proposition}}
\newcommand{\bt}{\begin{theorem}}
\newcommand{\et}{\end{theorem}}
\newcommand{\bex}{\begin{Example}}
\newcommand{\eex}{\end{Example}}
\newcommand{\bc}{\begin{corollary}}
\newcommand{\ec}{\end{corollary}}
\newcommand{\bl}{\begin{lemma}}
\newcommand{\el}{\end{lemma}}

\begin{document}

\title
[Non-uniqueness of convex bodies]
{Non-uniqueness of convex bodies with prescribed  volumes of  sections and projections}

\author{Fedor Nazarov, Dmitry Ryabogin and Artem Zvavitch}
\address{Department of Mathematics, Kent State University,
Kent, OH 44242, USA} \email{nazarov@math.kent.edu}\email{ryabogin@math.kent.edu}\email{zvavitch@math.kent.edu}

\thanks{The first author is supported in
part by U.S.~National Science Foundation Grant DMS-0800243. The second and third named authors are supported in part by U.S.~National Science Foundation Grant DMS-1101636. }


\begin{abstract}
We show that  if $d\ge 4$ is even, then one can find two essentially
different convex bodies such that the volumes of their maximal sections, central sections, and projections coincide for all directions.
\end{abstract}

\maketitle

\section{Introduction}

As usual,  a  {\it convex body} $K\subset {\mathbb R}^d$ is a compact convex subset  of ${\mathbb R}^d$ with non-empty interior.  We assume that $0\in K$. We consider   the {\it central section function} $A_K$:
\begin{equation}\label{A}
A_K(u)=\mathrm{vol}_{d-1}(K\cap u^{\perp}),\qquad u\in {\mathbb S}^{d-1},
\end{equation}
the {\it maximal section function} $M_K$:
\begin{equation}\label{m}
M_K(u)=\max_{t\in {\mathbb R}}\mathrm{vol}_{d-1}(K\cap(u^{\perp}+tu)),\qquad u\in {\mathbb S}^{d-1},
\end{equation}
and the {\it projection function} $P_K$:
\begin{equation}\label{p}
P_K(u)=\mathrm{vol}_{d-1}(K|u^{\perp}),\qquad u\in {\mathbb S}^{d-1}.
\end{equation}
 Here  $u^{\perp}$ stands for the hyperplane passing through the origin and orthogonal to the unit vector $u$, $K\cap (u^{\perp}+tu)$ is the section of $K$ by the affine hyperplane $u^{\perp}+tu$,  and
$K|u^{\perp}$ is the projection of $K$ to $u^{\perp}$. Observe that $A_K\le M_K\le P_K$.
It is well known, \cite{Ga},  that for origin-symmetric bodies  {\it each} of the  functions $M_K=A_K$ and  $P_K$   determines  the convex body $K\subset {\mathbb R}^d$ uniquely.  More precisely, either of the conditions
$$
M_{K_1}(u)=M_{K_2}(u)\qquad \forall u\in {\mathbb S}^{d-1},
$$
and
$$
P_{K_1}(u)=P_{K_2}(u)\qquad \forall u\in {\mathbb S}^{d-1},
$$
implies $K_1=K_2$, provided $K_1$, $K_2$ are origin-symmetric and convex.

In this paper, we address the (im)possibility of analogous results for not necessarily symmetric convex bodies.

It is well known, \cite{BF},  that on the plane there are convex bodies $K$ that are {\it not} Euclidean discs, but nevertheless  satisfy $M_K(u)=P_K(u)=1$  for all $u\in {\mathbb S}^1$.  These are the {\it bodies of  constant width} $1$.

In 1929 T. Bonnesen asked whether {\it every}  convex body  $K\subset {\mathbb R}^3$  is uniquely defined by $P_K$ and $M_K$, (see \cite{BF}, page 51).  We note that
in any dimension $d\ge 3$, it is not even known whether the conditions
$M_K\equiv c_1$, $P_K\equiv c_2$ are incompatible for $c_1<c_2$.

In  1969 V. Klee asked whether the condition $M_{K_1}\equiv M_{K_2}$ implies $K_1=K_2$, or, at least, whether the condition $M_K\equiv c$ implies that  $K$ is a Euclidean ball, see \cite{Kl1}.

Recently,  R. Gardner and V. Yaskin, together with  the second and the third named authors constructed two bodies of revolution $K_1$, $K_2$ such that $K_1$ is origin-symmetric, $K_2$ is not origin-symmetric, but $M_{K_1}\equiv M_{K_2}$, thus answering the first version of Klee's question but not the second one (see \cite{GRYZ}).

The main results we will present in this paper are the following.

\bt\label{klee}
If  $d= 4$, there exists a convex body of revolution $K\subset {\mathbb R}^d$ satisfying $M_K\equiv const$ that is not a Euclidean ball.
\et

\bt\label{Bon}
If  $d\ge 4$ is  even,  there exist two essentially different convex bodies of revolution $K_1$, $K_2\subset {\mathbb R}^d$ such that $A_{K_1}\equiv A_{K_2}$, $M_{K_1}\equiv M_{K_2}$, and $P_{K_1}\equiv P_{K_2}$.
 \et

Theorem \ref{klee} answers the question of Klee in $\R^4$, and Theorem \ref{Bon} answers the analogue of the question of Bonnesen in even dimensions.

\begin{remark}
Theorem \ref{klee} is actually true in all dimensions, but the construction for $d\neq 4$ is long and rather technical, so we will present it 
 in a  separate paper.
\end{remark}

We borrowed the general idea of the construction of the bodies $K_1$ and $K_2$ in Theorem \ref{Bon} from \cite{RY}, which attributes it to \cite{GV} and \cite{GSW}.
It 
can be easily understood from the following illustration.

\begin{figure}[ht]
\includegraphics[width=300pt]{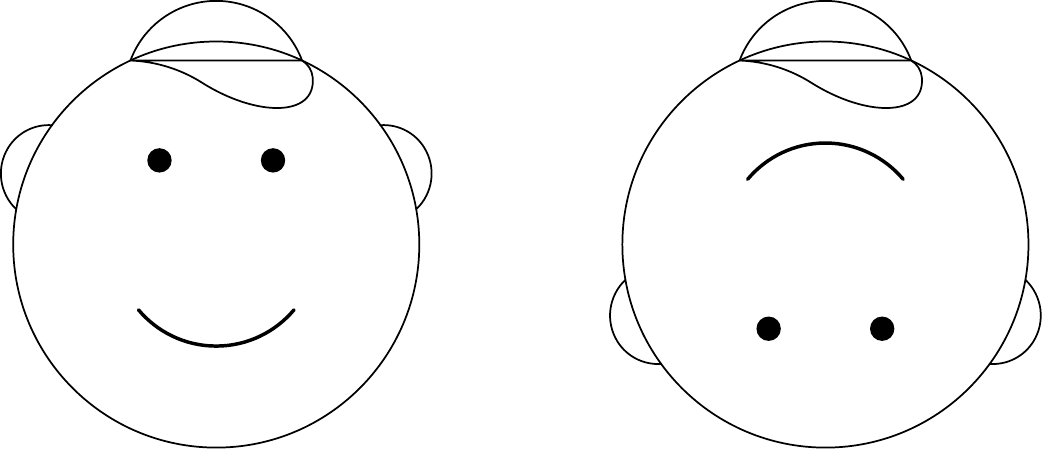}
\caption{Two small-eared round faces in a  cap}
\label{pic00}
\end{figure}

Here the "ears" and the "cap" will be made very small in order not to destroy the convexity of the bodies.

The paper is organized as follows. In  Section 2 we  reduce the problem to finding a non-trivial solution to two integral equations. In Section 3 we  prove Theorem \ref{klee}. In Section 4  we prove Theorem \ref{Bon}.

\section{Reduction to a system of integral equations}

From now on,  we assume that $d\ge 3$. We will be dealing with the   bodies of revolution
$$
K_f=\{x\in \R^d:\,x_2^2+x_3^2+...+x_d^2\le f^2(x_1)\},
$$
obtained by the rotation of   a smooth concave function $f$ supported on $[-1,1]$ about the $x_1$-axis.

Note that $K$ is rotation invariant, thus every its hyperplane section is equivalent to a section by a hyperplane with normal vector in the second quadrant of the $(x_1,x_2)$-plane.

\begin{figure}[ht]
\includegraphics[width=300pt]{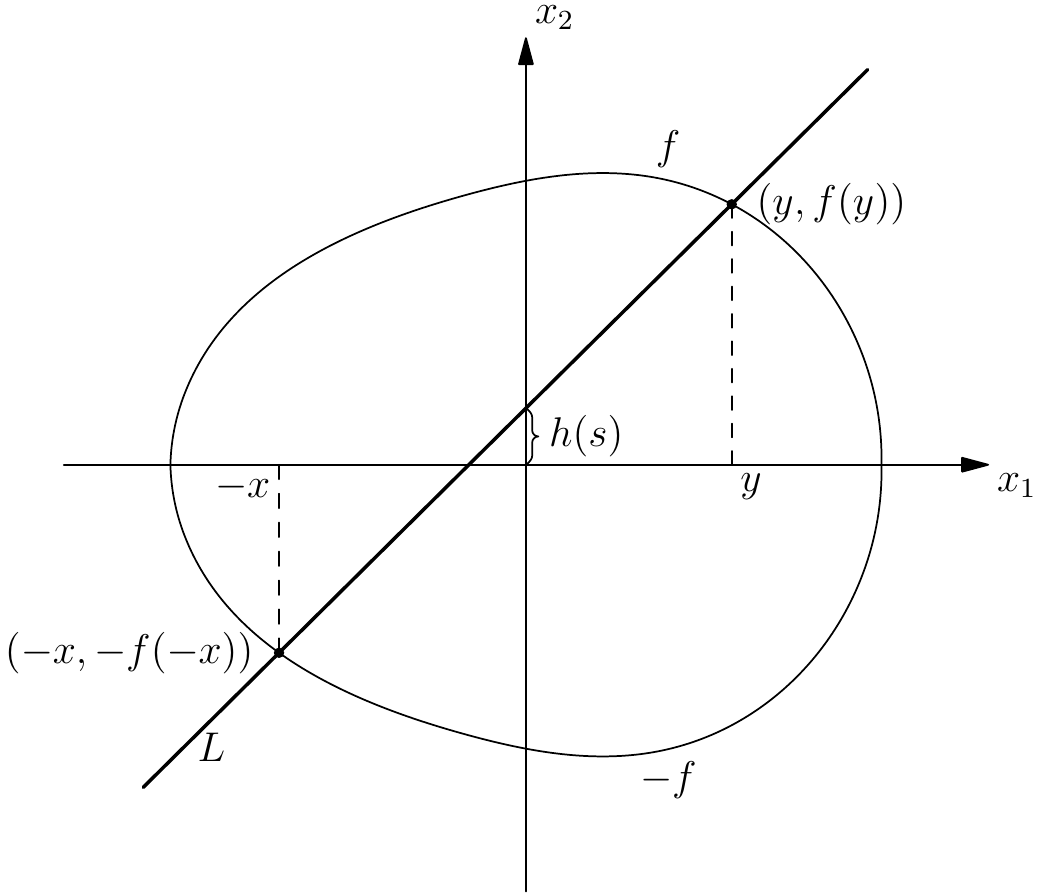}
\caption{View of $K$ and $H(L)$ in $(x_1,x_2)$-plane.}
\label{pic0}
\end{figure}

\bl\label{fn}
 Let $L(\xi)=L(s,h,\xi)=s\xi+h$ be a linear function with slope $s$, and let
  $H(L)=\{x\in{\mathbb R}^d:\,x_2=L(x_1)\}$ be the corresponding hyperplane.
  Then the   section $K\cap H(L)$ is of maximal volume iff
  \begin{equation}\label{max}
\int\limits_{-x}^{y}(f^2-L^2)^{(d-4)/2}L=0,
\end{equation}
where $-x$ and $y$ are the first coordinates of the points at which $L$ intersects the graphs of $-f$ and $f$ 
respectively (see Figure \ref{pic0}).
\el

\bp
Fix $s>0$. Observe  that
 the  slice $K\cap H(L)\cap H_{\xi}$ of $K\cap H(L)$ by the hyperplane $H_{\xi}=\{x\in {\mathbb R}^d:\,x_1=\xi\}$, $-1<\xi<1$,    is the
$(d-2)$-dimensional Euclidean ball $\{(x_3,x_4,...,x_d):\,x_3^2+...+x^2_d\le r^2\}$ of radius $r=\sqrt{f^2(\xi)-L^2(\xi)}$.
Hence,
\begin{equation}\label{const1}
\mathrm{vol}_{d-1}(K \cap H(L))=v_{d-2}\sqrt{1+s^2} \int\limits_{-x(s)}^{y(s)}(f^2(\xi)-L^2(\xi))^{(d-2)/2}d\xi,
\end{equation}
where $v_{d-2}$ is the volume of the unit ball in $\R^{d-2}$. 

The section $K\cap H(L)$ is  of maximal volume iff
$$
\frac{d}{d h}\mathrm{vol}_{d-1}(K \cap H(L))=0,
$$
where in  the {\it only if} part we use the Theorem of  Brunn, \cite{Ga}. Computing the derivative, we conclude that for a  given $s\in {\mathbb R}$, the section $K\cap H(L)$ is  of maximal volume iff (\ref{max}) holds.
\ep
\bl\label{fn1}
 Let $L(s,\xi)=s\xi+h(s)$ be a family of linear functions parameterized by the slope $s$. For each $L$ in our family, define the hyperplane $H(L)$ by
  $H(L)=\{x\in{\mathbb R}^d:\,x_2=L(x_1)\}$, (see Figure \ref{pic0}). The corresponding family of
 sections is of constant $d-1$-dimensional volume  iff
 \begin{equation}\label{const}
\int\limits_{-x}^{y}(f^2-L^2)^{(d-2)/2}=\frac{const}{\sqrt{1+s^2}}, \qquad \textrm{for all}\qquad s>0.
\end{equation}
\el
\bp
The right hand side in (\ref{const1}) is constant iff (\ref{const}) holds.
\ep

\section{The case $d=4$}

 Observe that when $d=4$, the system of equations (\ref{max}), (\ref{const}) simplifies to
  $$
  \int\limits_{-x}^yL=0,\qquad\text{and}\qquad \int\limits_{-x}^{y}(f^2-L^2)=\frac{const}{\sqrt{1+s^2}}, \qquad \textrm{for all}\qquad s>0.
  $$
In this case we will show that the {\it maximal sections correspond to level intervals}, see Proposition  \ref{prop1} below. We will also prove that {\it the values of the maximal section function $M_K$ depend on the distribution function $t\to |\{f>t\}|$ only}.  More precisely, we have

  \bt\label{distr}
 Let $d=4$,
$
K=\{x\in {\mathbb R}^4:\,x_2^2+x_3^2+x_4^2\le f^2(x_1)\}$, and let
$$
u=u(s)=(-\frac{s}{\sqrt{1+s^2}},\frac{1}{\sqrt{1+s^2}},0,0)\in {\mathbb S}^{3},\qquad s>0.
$$
 Then,
\begin{equation}\label{mkdistr}
M_{K_f}(u)=\pi\sqrt{1+s^2}\Big(\frac{2}{3}\,t^2\,|\{f>t\}|+\int\limits_t^{\infty}2\tau|\{f>\tau\}|d\tau\Big),
\end{equation}
where
$t$ is the unique solution of the equation $s=2t/ |\{f>t\}|$.

 In particular, if $f_1$ and $f_2$ are equimeasurable (i.e., for every $\tau >0$, we have $|\{ f_1>\tau \}|=|\{f_2>\tau \}|$), then $M_{K_{f_1}}\equiv M_{K_{f_2}}$.
 \et

\begin{figure}[ht]
\includegraphics[width=300pt]{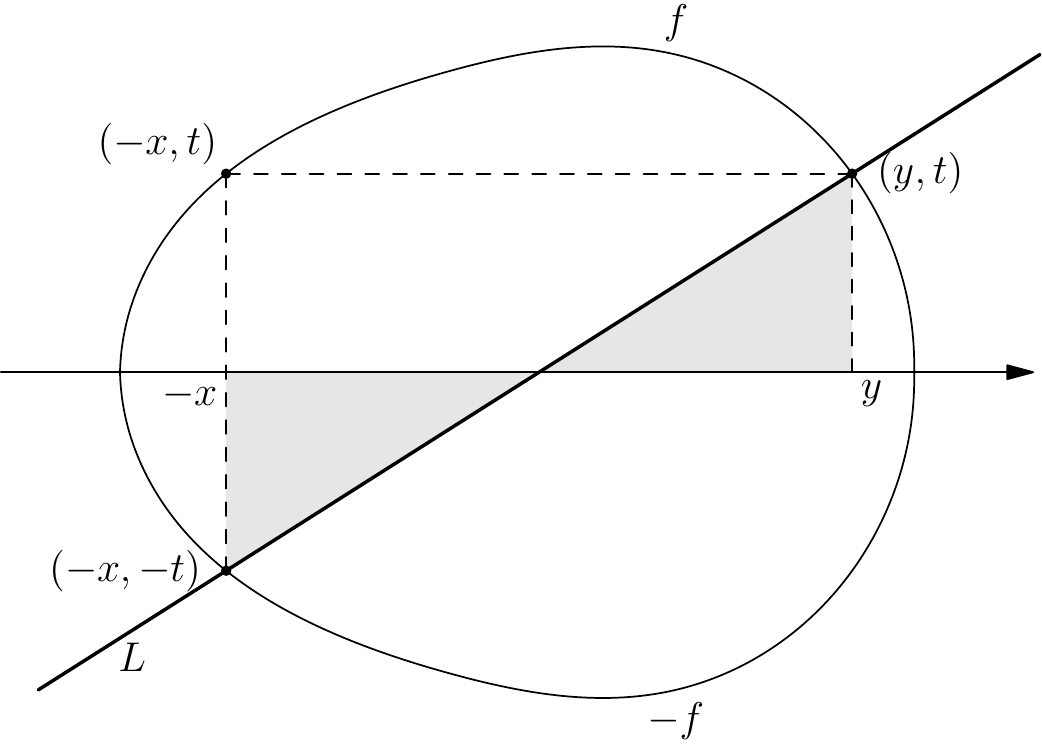}
\caption{Maximal slice in ${\mathbb R}^4$}
\label{pic1}
\end{figure}

Theorem \ref{distr} is a simple consequence of the following two propositions.
\bprop\label{prop1}
 Let $f, s, u(s)$ be as in Theorem \ref{distr}.   Then the section of  maximal volume
  in the direction $u(s)$  is the one that corresponds to the line joining $(-x,-t)$ and $(y,t)$, where $t$ is such that
 $s=2t/|\{f>t\}|$, $0< t< \max\limits_{\xi\in[-1,1]}f(\xi)$, (see Figure \ref{pic1}).
\eprop

\bp
Fix $s>0$.  Since the distribution function is decreasing to $0$,  there exists a unique $t$ satisfying $s=2t/|\{f>t\}|$. To prove that
$$
\int\limits_{-x}^yL(\xi)d\xi=0
$$
observe that two shaded triangles  on  Figure \ref{pic1} are congruent.
\ep
\bprop\label{prop2}
Let $K, f, t, x, y$ be as  in the previous proposition, and let the line $L$ be passing through the points $(-x,-t)$, $(y,t)$.
 Then (\ref{mkdistr}) holds.
\eprop
\bp
Note that
$$
\int\limits_{-x}^yL^2=(x+y)\frac{t^2}{3}=\frac{t^2}{3}|\{f>t\}|,
$$
and
$$
\int\limits_{-x}^yf^2=\int\limits_{\{f>t\}}f^2=t^2|\{f>t\}|+\int\limits_t^{\infty}2\tau|\{f>\tau\}|d\tau.
$$
\ep

{\bf Proof of Theorem 1}. Let $f_o(\xi)=\sqrt{1-\xi^2}$, $\xi\in[-1,1]$.
Take a concave function $f$ on $[-1,1]$ such that $f\neq f_o$ and $f$  is   equimeasurable with $f_o$.

\section{Proof of Theorem 2}

\begin{figure}[ht]
\includegraphics[width=350pt]{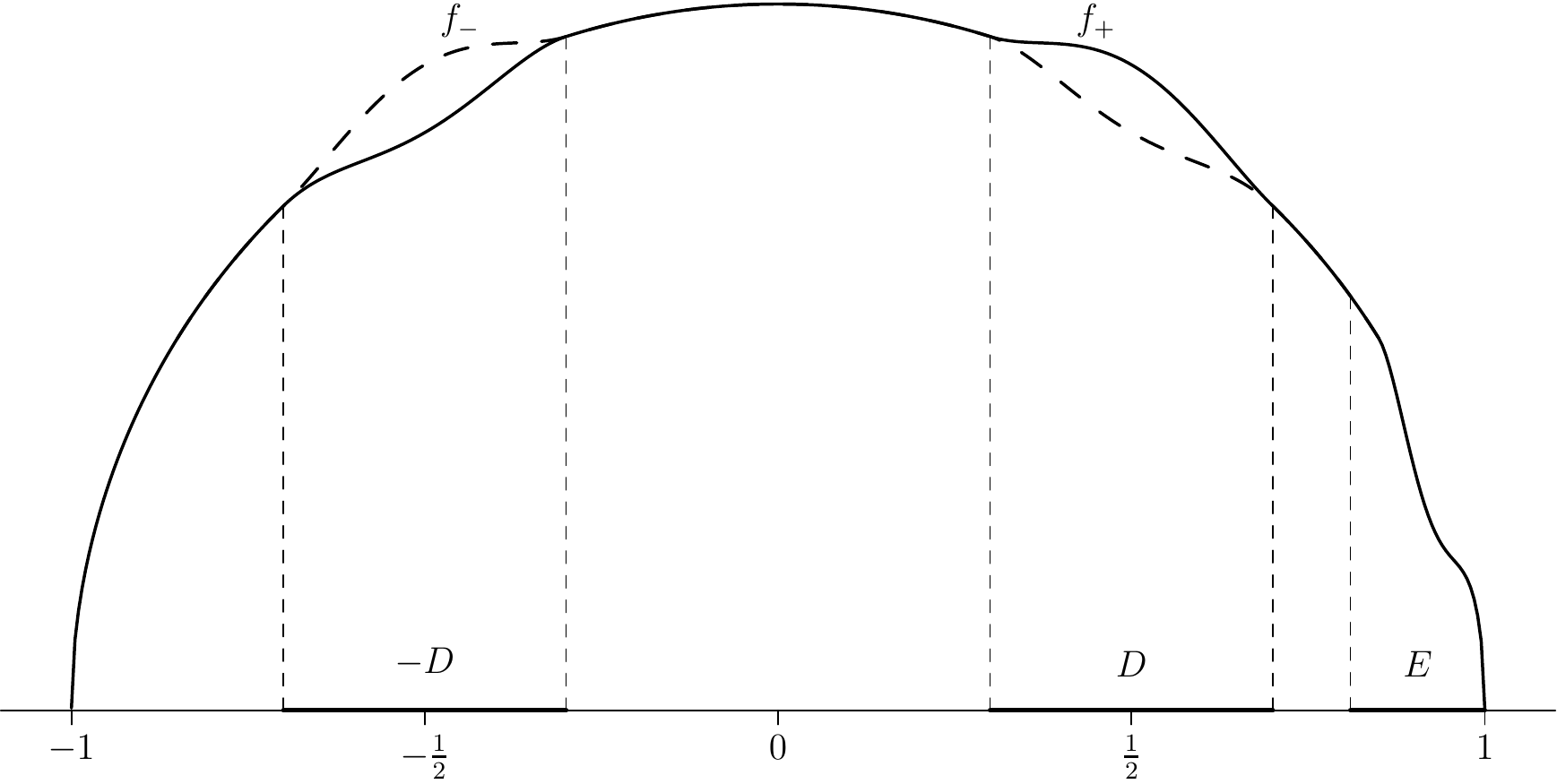}
\caption{Graph of functions $f_{\pm}$}
\label{pic3}
\end{figure}

Let $\varphi$ and $\psi$ be two  smooth   functions     supported on the intervals $D=[\frac{1}{2}-\delta, \frac{1}{2}+\delta]$ and  $E=[1- \delta,1]$ respectively,
where $0<\delta<\frac{1}{8}$.
Define
$$
f_+(\xi)=f_o(\xi)+ \varepsilon \varphi(\xi)-\varepsilon \varphi(-\xi)+\varepsilon \psi(\xi),
$$
and
$$
f_-(\xi)=f_o(\xi)- \varepsilon \varphi(\xi)+\varepsilon \varphi(-\xi)+\varepsilon \psi(\xi),
$$
where
  $\varepsilon>0$ is so small that $f_{\pm}$ are concave on $[-1,1]$  (see Figure \ref{pic3}).

Define $K_1=K_{f_+}$ and $K_2=K_{f_-}$.

Observe that
\begin{equation}\label{propnew}
f_+(\xi)=f_-(\xi)\,\,\forall \xi\in [-1,1]\setminus (D\cup(-D)),\,\,\textrm{and}\,\, f_+(\xi)=f_-(-\xi)\,\,\forall \xi\in D\cup(-D).
\end{equation}

We can choose $\varepsilon$ so small that $K_1$ and $K_2$ are very close to the Euclidean ball, and  the sections of the  maximal volume of $K_1$ and $K_2$
 are very close to the central sections of the ball.

In particular, the intersection points $x=x(s)$ and $y=y(s)$ satisfy
\begin{equation}\label{pojmal}
x(s),\,y(s)>\frac{5}{8}\,\,\, \text{if}\,\,\, s\le\frac{\sqrt{7}}{3},\qquad\text{ and}\qquad
x(s), \,y(s)<\frac{7}{8}\,\,\,\text{if}\,\,\, s\ge\frac{\sqrt{7}}{3},
\end{equation}
for both bodies.

First, we  show that $P_{K_1}\equiv P_{K_2}$ and $A_{K_1}\equiv A_{K_2}$.

Observe that we have $h_{K_1}(u)=h_{K_2}(u)$ and $\rho_{K_1}(u)=\rho_{K_2}(u)$ for all directions
$ u=(\xi,\sqrt{1-\xi^2},0,0,...,0)\in {\mathbb S}^{d-1}$,
 $\xi\in [0,1]\setminus D$. Observe also that  $h_{K_1}(u)=h_{K_2}(-u)$ and $\rho_{K_1}(u)=\rho_{K_2}(-u)$ for all directions
$ u=(\xi,\sqrt{1-\xi^2},0,0,...,0)\in {\mathbb S}^{d-1}$,
 $\xi\in  D$. Hence,
 the non-ordered  pairs
$$ \{h_{K_1}(u), h_{K_1}(-u)\},\qquad\text{and}\qquad \{h_{K_2}(u), h_{K_2}(-u)\}
$$
coincide for all $u\in{\mathbb S}^{d-1}$, and so do the pairs
$$
\{\rho_{K_1}(u), \rho_{K_1}(-u)\},\qquad\text{and}\qquad\{\rho_{K_2}(u), \rho_{K_2}(-u)\}.
$$

By the result of Goodey, Schneider and Weil,  \cite{GSW}, we have
  $P_{K_1}\equiv P_{K_2}$.
Also,
 $$
A_{K_1}(\theta)=\frac{1}{d-1}\int\limits_{{\mathbb S}^{d-1}\cap \theta^{\perp}} \frac{\rho_{K_1}^{d-1}(-u)+\rho_{K_1}^{d-1}(u)}{2}d\sigma(u)=
$$
$$
\frac{1}{d-1}\int\limits_{{\mathbb S}^{d-1}\cap \theta^{\perp}} \frac{\rho_{K_2}^{d-1}(-u)+\rho_{K_2}^{d-1}(u)}{2}d\sigma(u)=A_{K_2}(\theta)
$$
for all $\theta\in{\mathbb S}^{d-1}$.

It remains to show that
\begin{equation}\label{mdauzh}
M_{K_1}\equiv M_{K_2}.
\end{equation}
Assume  first that $d=4$.
The functions  $f_{+}$ and $f_-$ are equimeasurable, so (\ref{mdauzh}) follows from Theorem 3.

Let now $d\ge 6$ be even. Note that  in this case, $p=\frac{d-2}{2}\in{\mathbb N}$.

We claim that we can choose $\varphi$ such that (\ref{mdauzh}) holds. This result will be a consequence of the following two propositions.

\bprop\label{propnew1} If $\varepsilon$ is small enough, then  for
$$
u=u(s)=(-\frac{s}{\sqrt{1+s^2}},\frac{1}{\sqrt{1+s^2}},0,0,...,0)\in {\mathbb S}^{d-1}
$$
we have
\begin{equation}\label{mkdistr10}
M_{K_1}(u)=M_{K_2}(u),
\end{equation}
provided $s\ge \frac{\sqrt{7}}{3}$.
\eprop
 \bp
By (\ref{pojmal}), if the sections $K_1\cap H(L_1)$, $K_2\cap H(L_2)$ are the sections of the maximal volume, corresponding to the same slope $s\ge \frac{\sqrt{7}}{3}$,
then, they are the sections
of two symmetric (to each other) bodies corresponding to $\psi=0$.
 Hence, (\ref{mkdistr10}) holds by symmetry.
  \ep
To formulate the second proposition we will need the following result, which is a consequence of the Borsuk-Ulam Theorem.
\bl\label{mom11}
There exists $\varphi\neq 0$  such that
\begin{equation}\label{redbon}
\int\limits_{-\frac{5}{8}}^{\frac{5}{8}} f_+^{2j}(\xi)\xi^ld\xi= \int\limits_{-\frac{5}{8}}^{\frac{5}{8}} f_-^{2j}(\xi)\xi^ld\xi
\end{equation}
 for all  $j=0,...,p$ and $l=0,...,2(p-j)$.
\el
\bp
We will choose $\varphi$  using the
 Borsuk-Ulam Theorem.

For  $j$ and $l$ as above,  consider the vector $\mathbf{a}=\mathbf{a}(\varphi)$  with coordinates
$$
a_{j,l}(\varphi)=\int\limits_{-\frac{5}{8}}^{\frac{5}{8}}f_+^{2j}(\xi)\xi^ld\xi-\int\limits_{-\frac{5}{8}}^{\frac{5}{8}}f_-^{2j}(\xi)\xi^ld\xi.
$$
We will view this vector as an element of $\R^{n(p)}$ with appropriately chosen $n(p)$.

For $x=(x_0,x_1,...,x_{n(p)})\in {\mathbb S}^{n(p)}$  define $\varphi_x=\sum\limits_{j=0}^{n(p)}x_j\varphi_j$, where $\varphi_j$ are smooth not identically zero functions with pairwise disjoint supports
contained in $D$,
and let
$$
\mathbf{B}:\,(x_0,x_1,...,x_{n(p)})\,\,\to \,\,\varphi_x\,\,\to\,\,\mathbf{a}(\varphi_x)
$$
be  our map from ${\mathbb S}^{n(p)}$ to ${\mathbb R}^{n(p)}$.
By the definition of $f_{\pm}$ the map $\mathbf{B}$ is {\it odd}. Hence,  by the Borsuk-Ulam Theorem, one  can choose $x$ (and hence $\varphi_x$)
 in such a way that $\mathbf{a}(\varphi_x)=\mathbf{0}$.
\ep

\bprop\label{propnew2}
Let $L=L(s,h,\xi)$ be as above and let $0\le s\le \frac{\sqrt{7}}{3}$. Then,
\begin{equation}\label{devenconNEW}
\int\limits_{-x}^{y}(f_-^2-L^2)^p=\int\limits_{-x}^{y}(f_+^2-L^2)^p.
\end{equation}
Moreover,
\begin{equation}\label{devenmaxNEW}
\int\limits_{-x}^{y}(f_+^2-L^2)^{p-1}L =0
\qquad
 \text{if and only if}\qquad
   \int\limits_{-x}^{y}(f_-^2-L^2)^{p-1}L =0.
\end{equation}
In particular,  (\ref{mkdistr10}) holds for $0\le s\le \frac{\sqrt{7}}{3}$.
\eprop
\bp
We  start with the proof of  (\ref{devenconNEW}).
We open parentheses
and observe that all we need to prove is
$$
\int\limits_{-x(s)}^{y(s)} f_+^{2j}(\xi)\xi^ld\xi= \int\limits_{-x(s)}^{y(s)} f_-^{2j}(\xi)\xi^ld\xi
$$
for all $0\le s\le \frac{\sqrt{7}}{3}$, and for all $j=0,...,p$ and $l=0,...,2(p-j)$.
By (\ref{pojmal}), this follows from (\ref{redbon})  since $f_+(\xi)=f_-(\xi)$ for $|\xi|\ge \frac{5}{8}$.

 Similarly, (\ref{redbon}) implies (\ref{devenmaxNEW})  for $j=0,...,p-1$ and  $l=0,...,2(p-1-j)+1$.
\ep

Thus, (\ref{mdauzh}) follows from Propositions \ref{propnew1}, \ref{propnew2}. This finishes the proof of the Theorem.

\end{document}